\documentclass[preprint,12pt,sort&compress]{elsarticle}
\usepackage[left=2.0cm,right=2.0cm,top=1cm,bottom=2cm,bindingoffset=0cm]{geometry}
\usepackage{eepic}
\usepackage{amsmath,amssymb,amsthm}
\usepackage{multicol}
\usepackage{subfig}
\usepackage{enumitem}

\usepackage[T2A]{fontenc}
\usepackage[english]{babel}

\usepackage{listings, xcolor}

\begin{document}

\begin{frontmatter}

 \title{Forecasting and stabilizing chaotic regimes in two macroeconomic models via artificial intelligence technologies and control methods}

\author[hse]{Tatyana Alexeeva}%\corref{cor}\,}
%\ead{Corresponding author: tatyanalexeeva@gmail.com}
\author[d] {Quoc Bao Diep\,}
\author[spbu,ipm]{Nikolay Kuznetsov\corref{cor}\,}
\ead{Corresponding author: nkuznetsov239@gmail.com}
%\author[spbu]{Timur Mokaev\,}
\author[vsb]{Ivan Zelinka}

 \address[hse]{St.~Petersburg School of Physics, Mathematics and Computer Science, HSE University, \\194100 St. Petersburg, Kantemirovskaya ul., 3, Russia}
 %\address[kan]{Department of Economics, University of Kansas, Lawrence, KS 66045, USA}
 \address[d]{Faculty of Mechanical-Electrical and Computer Engineering, School of Technology, \\Van Lang University, Ho Chi Minh City, Vietnam}
 %\address[cefs]{Center for Financial Stability, New York, NY 10036, USA}
 \address[spbu]{Faculty of Mathematics and Mechanics,
 St. Petersburg State University, 198504 Peterhof,
 St. Petersburg, Russia}
 %\address[fin]{Department of Mathematical Information Technology,
 %University of Jyv\"{a}skyl\"{a},  40014 Jyv\"{a}skyl\"{a}, Finland}
 \address[ipm]{Institute for Problems in Mechanical Engineering RAS, 199178 St. Petersburg, V.O., Bolshoj pr., 61, Russia}
 \address[vsb]{Faculty of Electrical Engineering and Computer Science, V\v{S}B-TUO,\\ 17.listopadu 2172/15, 708 00 Ostrava-Poruba, Czech Republic}

\begin{abstract}
% The study of the complex behavior of economic processes associated with nonlinearity and irregular dynamics, arising as a result of internal mechanisms or external shocks in the mathematical models describing them, is an important task in current economic research. Irregular behavior can reduce the predictive properties of the model.
% In our article, we demonstrate how to overcome the difficulties associated with irregular, including chaotic, dynamics applying the evolutionary algorithms, the continuous deep Q-learning method and the delay feedback control method proposed by K.~Pyragas using two economic models as examples.
One of the key tasks in the economy is forecasting the economic agents' expectations of the future values of economic variables using mathematical models. The behavior of mathematical models can be irregular, including chaotic, which reduces their predictive power. In this paper, we study the regimes of behavior of two economic models and identify irregular dynamics in them. 
Using these models as an example, we demonstrate the effectiveness of evolutionary algorithms and the continuous deep Q-learning method in combination with Pyragas control method for deriving a control action that stabilizes unstable periodic trajectories and suppresses chaotic dynamics. We compare qualitative and quantitative characteristics of the model's dynamics before and after applying control and verify the obtained results by numerical simulation.
Proposed approach can improve the reliability of forecasting and tuning of the economic mechanism to achieve maximum decision-making efficiency.
\end{abstract}


\end{frontmatter}

\section{Introduction}

The processes taking place in the modern economy are characterized by increasing complexity and unpredictability, due to transformation of its structure, development and implementation of new technologies, as well as climate change and the epidemiological environment. Such challenges motivated posing a set of important problems in economic research as well as in actual economics activities:
   studying the complex behavior of models, for example, associated with non-linearity, large-scaling and a large number of variables or parameters (see, e.g. \cite{MaliarMW-JME-2021,FernandezNSLV-NBER-DL2020});
   forecasting irregular dynamics of economic processes, including critical transitions and crises (see, e.g. \cite{KindlebergerA-book2015});
   solving complex dynamic nonlinear models with a large number of heterogeneous agents in the presence of numerous stochastic shocks, for example Dynamic Stochastic General Equilibrium (DSGE) models (see, e.g. \cite{LindeSW-HB-2016,SlobodyanW-AEJM-2012}) 
   which are used by central banks to analyze and produce real or near-real time forecasting of important macro variables such as inflation and economic activity;
   monitoring of the degree of financial system's interconnectedness and stability, etc. \cite{Consoli-2021-BookDataSciEc};
   detecting patterns of behavior by different economic agents that are causing concerns (fraud, money laundering, collusion in procurement auctions, intentional default, etc.).

In economic research the tasks are focused on mathematical modeling of economic processes, forecasting the qualitative properties of their dynamics and controlling various regimes that arise in the economic mechanism, including under conditions leading to crisis states. In nonlinear mathematical models, signals of crisis conditions can manifest themselves in the form of irregular dynamics, including quasiperiodic and even chaotic ones. Additional complexity of the dynamics can be also associated with various unstable periodic trajectories embedded into the chaotic attractor of the dynamical system. Irregular behavior reduces the reliability of forecasts and thus undermines the predictive power of the model. In this case, decision-makers do not have the ability to predict and regulate the expectations of agents as well as evaluate the implications of economic policies.

Successfully overcoming these difficulties crucially depends on the choice of methods and effective ways of their application to the analysis and control of model dynamics. Artificial Intelligence (AI) technologies such as evolutionary algorithms (EAs) \cite{StornP-1997,davendrazelinka2016} and the reinforcement learning (RL) (the continuous deep Q-learning method and the actor-critic method) \cite{WatkinsD-QL-1992,SuttonB-RL-book1998,IkemotoU-QL-2019,IkemotoU-QL-2021} have significant advantages as optimization and adata-based control policy algorithms. However, we found that to correctly solve the problem of controlling irregular dynamics, including chaotic, and to reliably forecast the behavior of models, the usage of these technologies as sole method is not always sufficient. 
In some cases, the application of time-delay feedback control method (TDFC) proposed by K.~Pyragas \cite{Pyragas1992} as comlement to EAs or Q-learning is more effective approach.
Such combination of the methods allows targeted utilization of the control procedure to stabilize unstable periodic orbits (UPOs) of the model taking into account a problem statement in a subject field (in this case, in economics), the structure of the mathematical model and the limitations of the TDFC computational method.

We demonstrate the efficiency of this approach for two discrete-time macroeconomic models as examples, presenting two possible options for combining one of the AI technologies and the Pyragas method. In the first case, we apply a combination of EAs with the Pyragas method to one of overlapping generations (OLG) models \cite{Diamond-1965,Samuelson-1958-JPE} that an ``internal'' control variable  is included as economic feature \cite{AlexeevaDKMZ2-AI-2022}.
Using the computational abilities of a supercomputer, we show how the combination of differential evolution (DE) \cite{price2013differential} and self-organized migration (SOMA) algorithms \cite{zelinka2016soma} with the Pyragas method can overcome the problems associated with fractional-power nonlinearities in economic variables, significantly increase the performance of chaotic dynamics control and allow us to obtain much faster and more accurate fine-tuning of control parameters to achieve the desired state of the model and improve its forecasting behavior.
Our second model describes spatio-temporal (ST) model of pricing patterns in the global network market of goods \cite{Kim-2009-henon} that can demonstrate a rich spectrum of nontrivial dynamics, including multistability and chaotic regimes, which stimulate the use of increasingly sophisticated analytical and numerical techniques and tools based on artificial intelligence technologies. Using this model as an example, we show how to model ``external'' control and stabilize the chaotic regime by the Q-learning procedure and the Pyragas method.

\section{Problem statement}

Irregular behavior is an undesirable phenomenon in the economy, since the model of economic phenomenon has some interpretative power for making decisions and affecting the agents' expectations only if model possesses determinate (unique) economic equilibria. Therefore, revealing regular regime of the model behavior and regularization of the dynamics in the case of an irregular regime using control is one of the important tasks of modeling \cite{Pyragas-2006,OttGY-1990,BoccalettiGLMM-2000}.
In this paper, we consider two economies represented by an OLG model and a spatio-temporal (ST) model of pricing patterns in the global network market of goods; both models can demonstrate irregular, including chaotic, dynamics  for some parameter values. To control chaos we may find and stabilize an UPO embedded within a chaotic attractor. To solve this problem we refer to the time-delay feedback control (TDFC) approach, suggested by K.~Pyragas.  
The main idea behind the Pyragas approach is to stabilize UPO by constructing a control force proportional to the difference between the current state of the system and an earlier state of the system (delayed by some time interval). 
In this paper, we represent mentioned above models in the following extended mathematical form:
\begin{equation}
v_{t+1} =\varphi\left(v_t\right) + u_t\!, \quad %v(0)=v_0,\, v_0\in \mathcal{V},\, 
t\in \mathbb{Z}_{+},
\label{gen-sys}
\end{equation}
where $v_t\in \mathcal{V} \subseteq \mathbb{R}^n$, $\varphi\left(v_t\right)$ is a vector-function, $u_t$ is control input.
In the spirit of the ideas by Pyragas, we can consider $u_t$ as the Pyragas-like control:
\begin{equation}
u_t=u\left(v_t,v_{t-1},\ldots,v_{t-m}\right),
\label{pyragas-like-control}
\end{equation}
 such that along UPO with the period-$m$  (i.e. for $v_t=v_{t-m}$, $m>0$) the following relation holds:
\begin{equation}
u\left(v_t,v_{t-1},\ldots,v_{t-m}\right)=0.
\label{control-u0}
\end{equation}
 For the original Pyragas method, we have:
\begin{equation}
u\left(v_t,v_{t-1},\ldots,v_{t-m}\right)=K\left(v_t-v_{t-m}\right)\!,
\label{pyragas-orig}
\end{equation}
where the feedback gain $K$ can be represented either a scalar or a constant matrix.

For the OLG model, control variable is naturally embedded in the model in the form of a proportional tax, thus, it is inherently ``internal'' economic characteristic. Therefore, we can stabilize the chaotic dynamics in the OLG model using a sufficiently small control by adjusting the proportional tax a little.
The pricing model in the global network market does not contain a control variable explicitly, so we introduce an additive ``external'' control variable into this model. Next, using mathematical modeling we synthesize control variables in the Pyragas form in such a way that the current unstable trajectory could be attracted to a periodic orbit, upon reaching which the control variable is set to 0. To effectively optimize and fine-tune the control parameters we apply evolutionary algorithms (EAs) in the OLG model and the continuous deep Q-learning and the actor-critic method in the pricing model.

 \section{Pyragas method and evolutionary algorithms: implementation in~OLG~model}

We consider a two period OLG model with production and endogenous labor choice. Similar OLG models allow to exploring important intergenerational mechanisms such as demographic trends, the dynamics of education and retirement, capital accumulation, public policies, inflation, fiscal policies, etc. Two types of agents, consumers and firms, solve their dynamical optimization problems.\\
$(1)$ {\it Consumer problem}
\begin{equation}
%\begin{align*}
\max \text{ }U_{t}   =U_{1}\left(  c_{t}^{t}\right)  +U_{2}\left(  c_{t+1}%
^{t}\right)  -\bar{U}\left(  l_{t}\right)  ,
%\end{align*}
\label{cp:obj-function}
\end{equation}
\begin{equation}
\begin{aligned}
%&s.t.\\
&\mbox{s.t.} \qquad &c_{t}^{t}+{\kappa}_{t}   =\left(  1-\rho_{t}\right)  w_{t}l_{t},\\
&&c_{t+1}^{t}  =\tilde{R}_{t+1}\kappa_{t},
\end{aligned}
\label{cp:constraints}
\end{equation}
where $c_{t}^{t}=c_{t}$ is consumption of young in period $t$,
$c_{t+1}^{t}=c_{t+1}$ is consumption of young in period $t+1$ when they
will be old; 
$U_{1}\!\left(c_{t}^{t}\right)$ is the utility of consuming while
young, $U_{2}\!\left(c_{t+1}^{t}\right)$ the (future) utility of consuming
when old, and $\bar{U}\left(l_{t}\right)$ disutility of labor. In the first
period, the {\it budget constraint} \eqref{cp:constraints} says that consumption $c_{t}^{t}$ plus capital
$\kappa_{t}$ equals the income of the young, given by their labor income
$w_{t}l_{t}$ net of proportional taxes with rate $\rho_{t}$.
Here  $l_t$ are working hours (labor) when young, $w_t$ is the real wage rate, $\tilde{R}_{t+1}$ is the gross real interest rate at the relevant time.
In the second period of their life (physical time $t+1$), the agents born at $t$ consume their savings with interest. Non-negativity constraints form part of the model.\\
$(2)$ {\it Firms' problem}
\begin{equation}
\begin{aligned}
\max \text{ }\underset{\kappa_{t},l_{t}}{\Pi_{t}} = l_{t}^{\beta} {\kappa}_{t-1}^{1-\beta} - w_t l_t - \tilde{R}_t {\kappa}_{t-1}
,
\label{fp:max_obj_func_profit}
\end{aligned}
\end{equation}
where $\Pi_{t}$ is the profit, $0 < \beta \leq 1$ is the parameter of the Leontieff production function.

The solution of this OLG model with respect to variables $(c_t,l_t)$ is described by the following two-dimensional map
$\psi : \mathbb{R}^2 \to \mathbb{R}^2$
\begin{equation}\label{eq:olg}
    \psi(c_t,l_t) = \left(l_t^\sigma - c_t, \ b_L ( \beta l_t - c_t^\lambda - g_t)\right),
\end{equation}
which generates the following discrete-time dynamical model 
\begin{equation}\label{eq:olg:lambda}
    \begin{cases}
        c_{t+1} = l_t^\sigma - c_t,\\
        l_{t+1} = b_L (\beta l_t - c_t^\lambda - g_t),
    \end{cases} \qquad
    t \in \mathbb{Z}_+,
\end{equation}
where 
$g_t$ is the government spending in period $t$;
 $b_L > 1$ is the parameter of the Leontieff production function, $\lambda >1$ and $\sigma \geq 1$ are model parameters (elastisities of functions).

Model \eqref{eq:olg:lambda} describes complex behavior of the agents in conditions of economic equilibrium -- a situation when supply and demand in all markets are balanced. 
Dynamics of agent's consumption and labor in the model can be irregular, involving chaotic regime, that decreases predictable power of the model and agents' expectations. In this case, one of the effective ways to suppress chaotic behavior in \eqref{eq:olg:lambda} through stabilizing UPOs embedded in a chaotic attractor using a minimal control.
The latter task requires to determine any period-$m$ UPO and introduce a control in the model.
To solve this problem and provide the minimum control through the variable $g_t$ in \eqref{eq:olg:lambda} we choose $g_t$ in the Pyragas form \eqref{pyragas-orig} ($g_t$ thus equals 0 on the target trajectory, i.e. UPO). Within framework of the economic problem, $g_t$  can be generated as a time series and implemented as follows:
\begin{equation}\label{eq:tax}
g_t=\beta \rho_t l_t,
\end{equation}
where $\rho_t$ is the rate of the proportional tax on the young.
  In the absence of a control variable (i.e. $g_t=0$), we get the following
\begin{equation}\label{eq:olg-out-control}
    \tilde{\psi}(c_t,l_t) = \left(l_t^\sigma - c_t, \ b_L ( \beta l_t - c_t^\lambda)\right).
\end{equation}
 
For system \eqref{eq:olg-out-control}, let us find a period-$m$ UPO such that $c_t=c_{t-m}$, $l_t=l_{t-m}$. Here, the original Pyragas approach \eqref{pyragas-like-control}--\eqref{pyragas-orig} does not allow us to identify and stabilize the UPO using only one control parameter. To overcome this difficulty we solve an auxiliary mathematical problem 
%within the framework of controlling economic parameters 
by adding third equation and artificial variable $\tilde{l}_t$ to system \eqref{eq:olg:lambda}, and consider the following extended system:  
\begin{equation}\label{eq:olg:dfc3}
    \begin{cases}
        c_{t+1} = l_t^\sigma - c_t,\\
        l_{t+1} = b_L \big(\beta l_{t} - c_{t}^\lambda\big) + b_Lk_1 \big(l_{t} - \tilde{l}_{t-(m-1)}\big),\\
        \tilde{l}_{t+1} = l_{t} + k_2 \big(l_{t} - \tilde{l}_{t-(m-1)}\big),        
    \end{cases}
\end{equation}
where $\tilde{l}_{t}$ specifies the delay, control parameters $k_1, k_2$ use only for stabilization UPO and pick up by TDFC method. That allows obtaining the control variable $g_t$ through economic parameters according to \eqref{eq:tax} and adjusting the behaviour of the economy using a fairly small disturbance inside model. 
Along the desired period-$m$ UPO in \eqref{eq:olg} or \eqref{eq:olg:lambda} the following condition holds $l_{t} - l_{t-m} = 0$. Plugging this equality into second and third equations of \eqref{eq:olg:dfc3}, we receive $\tilde{l}_{t+1}=l_t=l_{t-m}$ and, consequently, $\tilde{l}_t=\tilde{l}_{t-m}$.
Therefore, to numerically analyze the dynamics of \eqref{eq:olg:dfc3} we choose the initial conditions $(c_0,\,l_0,\,\tilde{l}_{0})$
so that if they get captured in the UPO, i.e. $c_0=c_{-m},\,l_0=l_{-m}=\tilde{l}_{-(m-1)},\,\tilde{l}_{0}=\tilde{l}_{-m}=l_{-1}$, then $u_t\equiv 0$. 
System \eqref{eq:olg:dfc3} is the notation of general mathematical model \eqref{pyragas-like-control}
 in terms of OLG model with{}:
\begin{equation}\label{gen-eqs}
  v_{t} = \left(
    \begin{array}{c}
      c_{t} \\
      l_{t} \\
      \tilde{l}_{t} \\
    \end{array}
  \right)\!,\,
  \varphi\left(v_{t}\right) = \left(
    \begin{array}{c}
      l_{t}^{\sigma}-c_{t} \\
      b_L\left(\beta l_{t} - c^{\lambda}_{t}\right) \\
      l_{t} \\
    \end{array}
  \right)\!,\,
   u_t = \left(
    \begin{array}{c}
      0 \\
      b_L k_1 \left(l_{t} - \tilde{l}_{t-(m-1)}\right) \\
      k_2\left(l_{t}-\tilde{l}_{t-(m-1)}\right)  \\
    \end{array}
  \right)\!.
  \end{equation}

In this work, to determine parameters of time-delay feedback control within Pyragas procedure,
we used three most powerful and commonly used evolutionary algorithms:
%\begin{enumerate}
  \label{itm:EAs:DE} differential evolution algorithm (DE/rand/1/bin)~\cite{StornP-1997};
  \label{itm:EAs:SOMA} self-organized migration algorithm (SOMA)~\cite{zelinka2016soma};
  \label{itm:EAs:SOMAT3A} SOMA with Team To Team Adaptive strategy (SOMA T3A)~\cite{Diep-2019,DiepZDS-2019}.
%\end{enumerate}
Using these EAs and solving numerically the corresponding equation
$(c_t, l_t) = \tilde{\psi}^m(c_t, \, l_t)$, $m = 1,2,\ldots$,
we obtain that for system~\eqref{eq:olg-out-control} with parameters
$\lambda = 3$, $\beta = 0.99$, $\sigma = 1.03$, $b_L = 1.54$
it is possible to reveal a candidate for a period-$m$ UPO, where $m=5$, to be further stabilized, i.e.
\begin{multline}\label{eq:upo-5}
    (0.49786228048149456, 0.9336025046020485) \to \\
    (0.433817924634687, 1.2333289032631394) \to \\
    (0.807294938247322, 1.7546020061352825) \to \\  
    (0.9771534043359568,  1.8648192605089056) \to \\
    (0.9228564283533333, 1.4062615940651981).
\end{multline}

\medskip

To fine-tune $k_1$, $k_2$, we again use the EAs,
finding the minimum value of the cost function
${\rm CF}(k_1, k_2)$ defined as a spectral radius of the fundamental matrix for the
corresponding map $\tilde{\psi}^m$,
% $\Phi(c_t,l_t) =  J(c_{t+4},l_{t+4}) \cdot J(c_{t+3},l_{t+3}) \cdot \ldots \cdot J(c_t,l_t)$
computed along the points of the UPO in \eqref{eq:upo-5} and reflecting its orbital stability
after stabilization.
% Here matrix $J(c_t,l_t)$ is the Jacobian of~\eqref{eq:olg:dfc}.
As a result, we obtain the following values: $k_1 = -0.18808852$, $k_2 = -0.95598289$.

\begin{figure}[!ht]
    \centering
    \includegraphics[width=\textwidth]{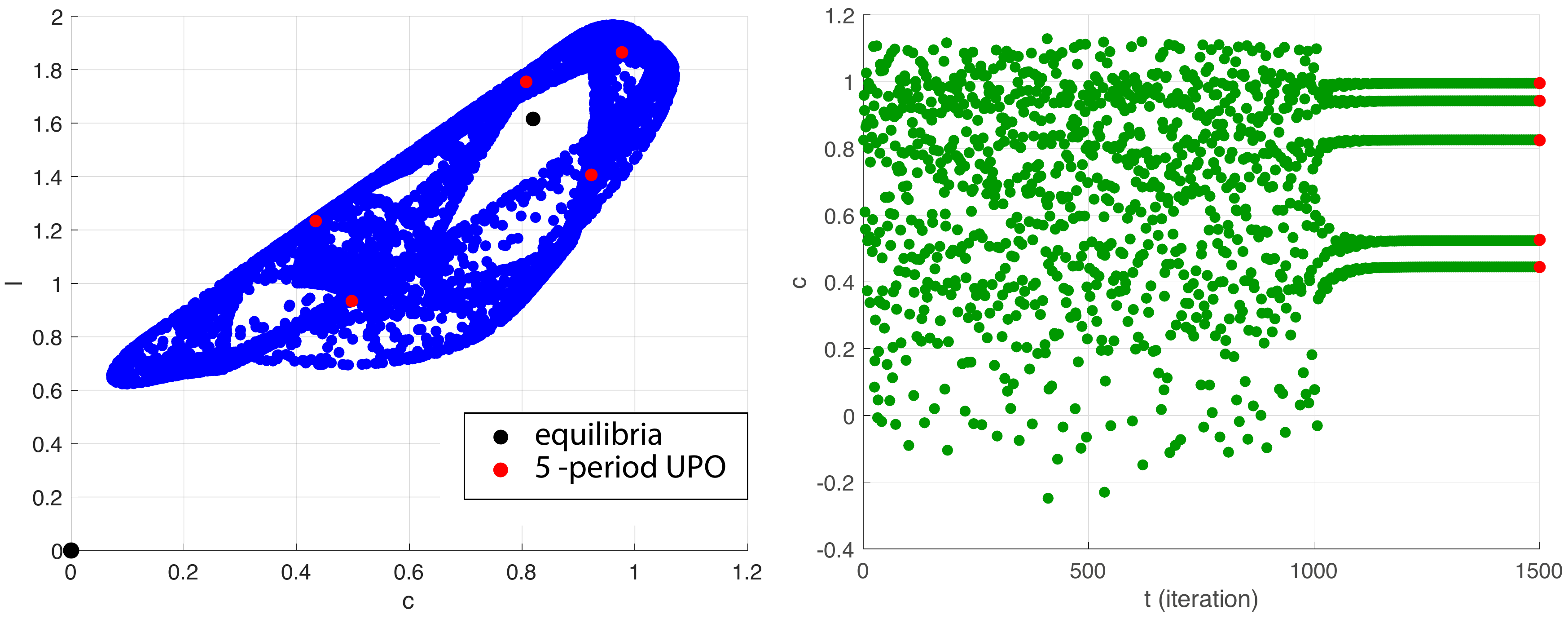}
    \caption{Period-5 UPO (red) embedded in the chaotic attractor (blue)
    in the OLG model \eqref{eq:olg-out-control} and it's stabilization (green) via the Pyragas control.}
\label{fig::olg:sc:chaos_upo5}
\end{figure}

Remark that the obtained values deliver stabilization only locally, i.e.
for some initial points in a vicinity of the UPO.
So, the next task in our experiments was to solve an \emph{optimal control problem} of
finding optimal values of control parameters $k_1$, $k_2$
to maximize the basin of attraction of the stabilized UPO.
To this end, we consider two regions:
for control parameters $(k_1, k_2) \in (-1, 1) \times (-1, 1)$, and
for system's parameters $(c, l) \, \in \, [0, 2] \times [0, 3]$
together with the specific partition step $c_{\rm step} = l_{\rm step} = \varepsilon$.
We define a specific grid of points
\[
  \mathcal{B}_{\rm grid}(\varepsilon) =
  \mathcal{B}_{\rm inf}(\varepsilon, k_1, k_2) \, \cup \,
  \mathcal{B}_{\rm upo}(\varepsilon, N_{\rm iter}, k_1, k_2) \, \cup \,
  \mathcal{B}_{\rm other}(\varepsilon, N_{\rm iter}, k_1, k_2),
\]
which consists of points, from which trajectories go to infinity
\big(set $\mathcal{B}_{\rm inf}(\varepsilon, k_1, k_2)$\big);
from which trajectories after $N_{\rm iter}$ iterations tend to the stabilized UPO
\big(set $\mathcal{B}_{\rm upo}(\varepsilon, N_{\rm iter}, k_1, k_2)$\big);
from which trajectories after $N_{\rm iter}$ iterations tend to some other attractors
\big(set $\mathcal{B}_{\rm other}(\varepsilon, N_{\rm iter}, k_1, k_2)$\big).
To use the computational abilities of EAs, we write this optimal control problem
in the following form with the corresponding cost function:
${\rm CF}(k_1, k_2) =
    \big|\mathcal{B}_{\rm grid}(\varepsilon)\big|
    \! - \! \big|\mathcal{B}_{\rm upo}(\varepsilon, N_{\rm iter}, k_1, k_2)\big| \ \to \ \min$.

Solving such a computational problem implies calculation
and examination of the behavior for the number of trajectories
of system~\eqref{eq:olg:dfc3} equal to $|\mathcal{B}_{\rm grid}(\varepsilon)|$
at each iteration of the evolutionary algorithm.
As $\varepsilon$ decreasing
and $N_{\rm iter}$ increasing, it turns into an extremely time and resource consuming
computational procedure.
In order to make this procedure more efficient, we implement it
in MathWorks Matlab using Parallel Computing Toolbox and launch it on two
powerful HPCs at IT4Innovations National Supercomputing Center of the Czech Republic, i.e.
Barbora Cluster and Salomon Cluster. 
%(for more details see~\cite{BarboraHPC-info} and~\cite{SalomonHPC-info}, respectively).
It took 324 Matlab workers (Matlab computational engines) in total with 48 hours of calculation time at each worker.

For stabilized UPO \eqref{eq:upo-5} and procedure parameters $\varepsilon = 0.01$, $N_{\rm iter} = 10000$
% (see Table~\ref{table:olg:upo5} and Fig.~\ref{fig::olg:basin:upo})
using three evolutionary
algorithms~\ref{itm:EAs:DE}, \ref{itm:EAs:SOMA}, \ref{itm:EAs:SOMAT3A}
we obtain that the best value of cost function
(i.e. the maximum number of points in
$\big|\mathcal{B}_{\rm upo}(\varepsilon, N_{\rm iter}, k_1, k_2)\big|$) equals $5728$ and
is reached at the point $(k_1, k_2) = (-0.1473, -0.8719)$.
This value is obtained by SOMA T3A algorithm launched on Salomon Cluster.
The simulation consumed about $150,000$ core hours for the computation, which is
equivalent to around $17.36$ continuously working years of the single CPU with $2.5$ GHz.

\begin{figure}[!ht]
\centering
\noindent\makebox[\textwidth]{%
\includegraphics[width=0.6\textwidth]{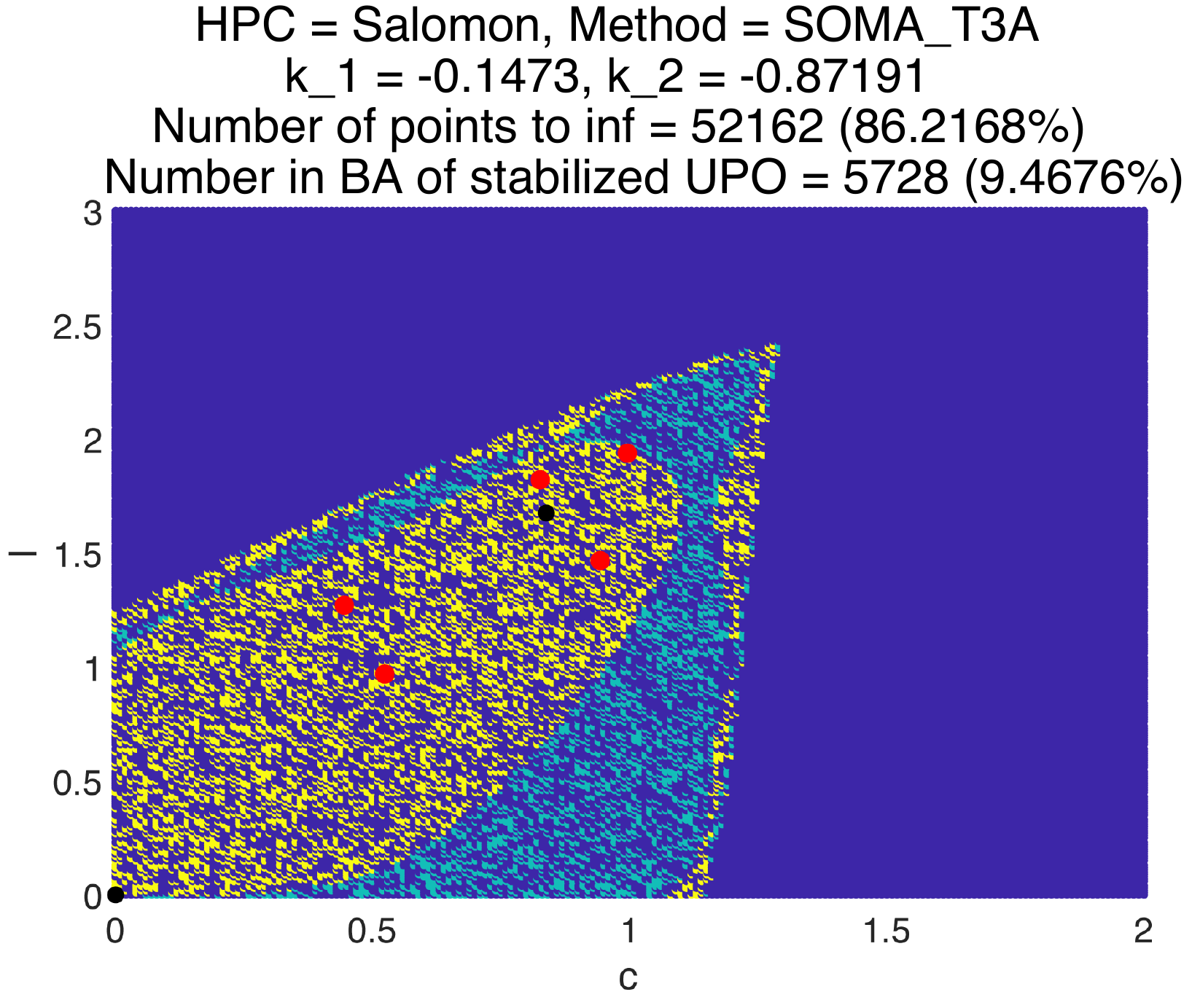}}
\caption{Basin of attraction (yellow) of stabilized UPO (red)
with respect to optimal control parameters $(k_1, k_2)$ for system~\eqref{eq:olg:dfc3}.
Purple domain corresponds to the points, from which trajectories go to infinity;
cyan domain corresponds the points, from which trajectories go to other attractors.}
\label{fig::olg:basin:upo}
\end{figure}

\section{Deep reinforcement learning: implementation in ST-model of pricing and H\'{e}non map}

One of the natural economic mechanisms in which irregular dynamics is formed is the pricing in the markets for goods that could only be stored for short-term, for example, in the markets for trade in fish and agricultural products. Such effects can be described by any type of cobweb model with non-linear supply and/or demand functions and various types of expectations of economic agents (naive, rational, adaptive). In addition, in the economic context, of interest may be the spatiotemporal behavior of the model associated with the interaction of local cobweb-type markets. Prices in such a model may demonstrate irregular dynamics, including chaotic ones. 

In our work, we represent an economy in which the global market consists of a network of local markets represented by two types of alternating nodes: ``even'' (or ``recipients'') and ``odd'' (or ``donors''), i.e. markets with a higher and a lower preference in the consumption of a good. This leads to a model of global trade in which the stationary net inflow of goods into different types of markets has opposite directions. 

We explain the global market dynamics by Lattice Dynamical System (LDS) \cite{Chow-LDS-2003,Kim-2009-henon}. 
At each node $j$ of a $d-$dimensional lattice $\mathbb{Z}^d$
we consider a finite dimensional local dynamical systems as map $f_j \colon M_j \to M_j$,
$M_j=\mathbb{R}^1,\, \forall j\in \mathbb{Z}$\, ($d=1$)
is a local phase space at the node $j$ (or an infinite chain), $f_j=f$.
Price dynamics could be formulated as follows:
\begin{equation}
    p_j^{(\tau+1)} = \left(\Phi\left(p^{(\tau)}\right)\right)_j
    = (1-\tilde\alpha) p_j^{(\tau)} + \tilde\alpha f\left(p_j^{(\tau)}\right)
    + \zeta \Bigl(p_{j-1}^{(\tau)} - 2 p_j^{(\tau)} + p_{j+1}^{(\tau)}\Bigr),
\label{LDS}
\end{equation}
where $p_j^{(\tau)}$ is good's price at node $j$
and at time $\tau>0$, $\left(\Phi\left(p^{(\tau)}\right)\right)_j=F\left(\{p_{j}^{(\tau)}\}^k\right)$
is the evolution operator, $\{p_{j}^{(\tau)}\}^k=\{p_i\mid |i - j|\leq k,\, k\geq 1\, \mbox{is integer}\}$,
$F:\mathbb{R}^{2k+1}\longrightarrow \mathbb{R}^k$ is a differentiable map
of class $C^2$, $\tilde{\alpha}\in (0,\,1]$ indexes the degree of nonlinearity in the local market dynamics
$f\left(p_j^{(\tau)}\right)$, $f=1-\Lambda p_{j-1}^{(\tau)}\left(1-p_{j-1}^{(\tau)}\right)$
is the logistic map with the parameter $\Lambda$, $\zeta>0$ is the diffusion coefficient measuring the intensity of interaction between the neighboring local markets.
A solution $p^{(\tau)}=\{p_j^{(\tau)}\}$ is stationary if $p_j^{(\tau)}=p_j=p_t$, $\forall t\in \mathbb{Z_{+}}$.
Then, the stationary solution of the global market dynamics could represent as the following
\begin{equation}
p_{t+1} = \frac{1}{\gamma_2}\Bigl(-\tilde\beta - \gamma_1 p_{t-1} + \left(2 + \tilde\beta\left(1+\Lambda\right)\right)p_t - \tilde\beta\Lambda p_t^2\Bigr),
\label{st-state-eq1}
\end{equation}
where $\tilde\beta=\frac{\tilde\alpha}{\zeta}$.
Taking $p_t\rightarrow d_1x_t+d_2,$\, $p_{t-1}\rightarrow d_1y_t+d_2$
$\left(d_1=\frac{\gamma_2}{\tilde\beta \Lambda}\right.,$\, $\left. d_2=\frac{2+\tilde\beta\left(1+\Lambda\right)}{2\tilde\beta\Lambda}\right)$ we can transform \eqref{st-state-eq1} to the 2D dynamical system
\begin{equation}\label{henon-sys}
    (x_{t}, y_{t}) = H^t_{a,b}(x_0,y_0) = \underbrace{H_{a,b} \circ \cdots \circ H_{a,b}}_{t \ \text{times}} (x_0,y_0), \quad
    t \in \mathbb{Z}_+, \quad (x_0,y_0) \in \mathbb{R}^2,
\end{equation}
generated by the H\'{e}non map $H_{a,b}(x,y)=\left(a + b y - x^2,\, x\right)$
with
$a=\frac{1}{{\gamma_2}^2}\Bigl(\gamma_2\left(1+\frac{\tilde\beta(1+\Lambda)}{2}\right)
\left(\frac{2}{\gamma_2} + b - 1\right) + {\tilde\beta}^2\left(\frac{1-\Lambda}{2}\right)^2 -1 \Bigr)$,
$b=-\frac{\gamma_1}{\gamma_2}=-\tilde\gamma$,
$0<\tilde\gamma<1$, and $\Lambda\in (0,\,4]$.
It has the unstable equilibria $O_{\pm}=\left(x_{\pm},\,y_{\pm}\right)$:
\[
    x_{\pm}=\frac{1}{2}\Bigl(b-1\pm\sqrt{\left(b-1\right)^2+4a}\Bigr), \quad
    y_{\pm}=\frac{1}{2}\Bigl(b-1\pm\sqrt{\left(b-1\right)^2+4a}\Bigr).
\]

Our goal of chaos control is to stabilize \eqref{henon-sys} at one of its UPOs. To achieve that, we can rewrite general system \eqref{gen-sys} in terms of H\'{e}non map \eqref{henon-sys} with{}\footnote{* is the transpose of a matrix}:
\begin{equation}\label{henon-0}
    v_t = (x_t,y_t)^{\ast}, \quad \varphi\left(v_t\right) = H_{a,b}(v_t), 
\end{equation}
and following the main ideas from \cite{IkemotoU-QL-2019,IkemotoU-QL-2021}, apply approach based on a RL algorithm called \emph{continuous deep Q-learning method} to synthesize a small control $u_t$ and use it to stabilize UPO with period $m=2$ in \eqref{henon-sys}.
In the framework of RL $u_t$ is considered as a controller action, subject to a certain control policy, in the process of interacting with a controlled environment, i.e. system \eqref{henon-0}, in order to transition the system from the current state to the next one in the best way. Thus, we can say the control problem is to find a control policy (function) which associates with each state $v_t$ in \eqref{henon-0} a control action $u_t$ such that the given goal is achieved in an optimal way indicating the quality of the control policy. As the criterion for the quality of control action taken $u_t$ in state $v_t$ a Q-function $Q(v_t,u_t)$ is used. It can be expressed by the sum of expected rewards $r_t$ for the transition of the system from the current state to the next one at time after taking an control action. We consider reward at time $t$ in the form of the following function:
\begin{equation}\label{r-reward}
r_t = R(v_t, u_t) = -\left(v_t-\varphi^2\left(v_t\right)\right)^* \, S_{1} \, \left(v_t-\varphi^2\left(v_t\right)\right) - 
    u_t^* \, S_{2} \, u_t, 
\end{equation}
where $S_{1,2}$ are some positive definite matrices.
Accordingly, we solve the following problem: calculate such a stabilizing control $u_t$ in which the Q-function takes the maximum value if control was started from a randomly chosen initial state. Extending the system, by analogy with the above case, we can look for a control in the Pyragas-like form \eqref{pyragas-like-control} $u_t=u(v_t,v_{t-2})$.

\begin{figure}[!ht]
    \centering
    \includegraphics[width=\textwidth]{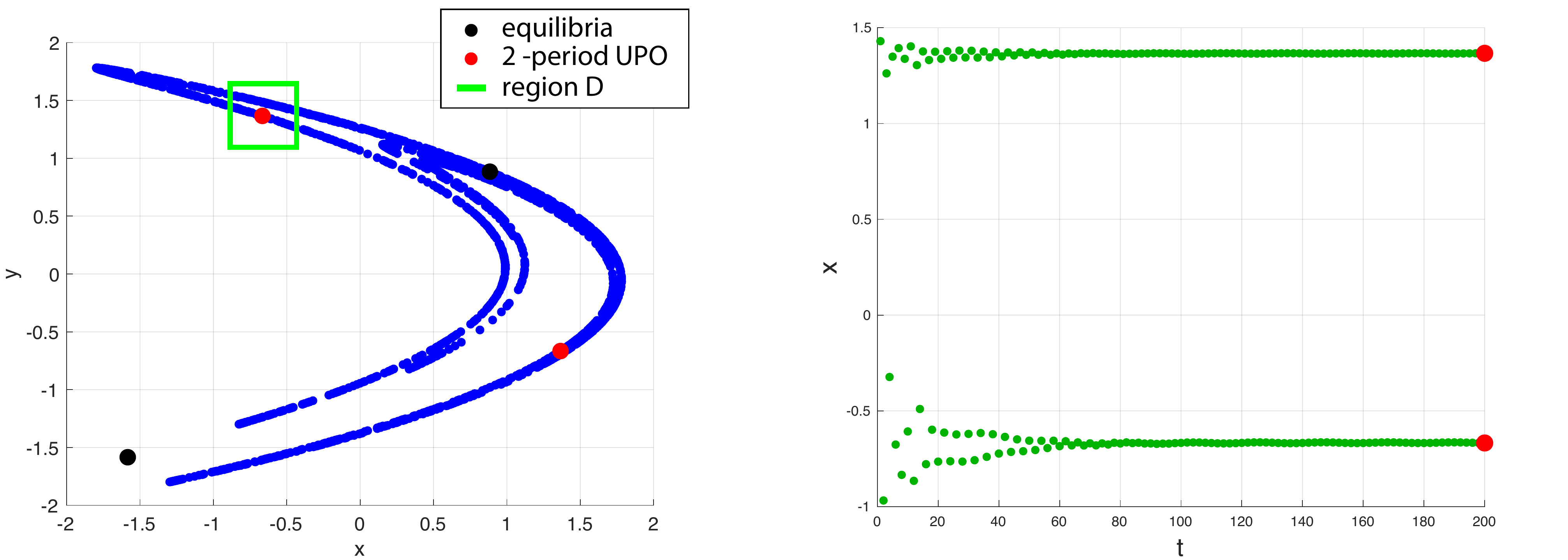}
    \caption{Period-2 UPO (red) embedded in the chaotic attractor (blue)
    in the H\'{e}non map \eqref{henon-sys} and it's stabilization (green).}
\label{fig::henon:chaos_upo2}
\end{figure}

The proposed approach based on the continuous deep Q-learning algorithm allows us to tune the optimal control more efficiently and stabilize the system faster compared to the direct application of control in the Pyragas form \eqref{pyragas-orig}, i.e. $u_t=K\left(v_t - v_{t-2}\right)$, where the selection of the control action $K$ depends significantly on the skills of the researcher.
For instance, for a H\'{e}non map \eqref{henon-sys} with ``classical'' parameter values $a = 1.4$, $b = 0.3$, we picked up the following matrix $K=\left[-1.1, 0; -0.9, 0.9\right]$ and considered the following initial conditions $v_0=[-0.6, 1.3]$, from which we were able to stabilize (see Fig.~\ref{fig::henon:chaos_upo2}) the following period-2 UPO:
\begin{multline*}
    (x^{\rm upo}_{1}, y^{\rm upo}_{1}) = (-.6661200716, 1.366120072) \to \\
    \to (x^{\rm upo}_{2}, y^{\rm upo}_{2}) =(1.366120072, -.6661200716)
\end{multline*}
in the region $D=[-0.85, -0.45]\times [1.1, 1.6]$.

\begin{figure}[!ht]
    \centering
    \noindent\makebox[\textwidth]{%
    \includegraphics[width=0.5\textwidth]{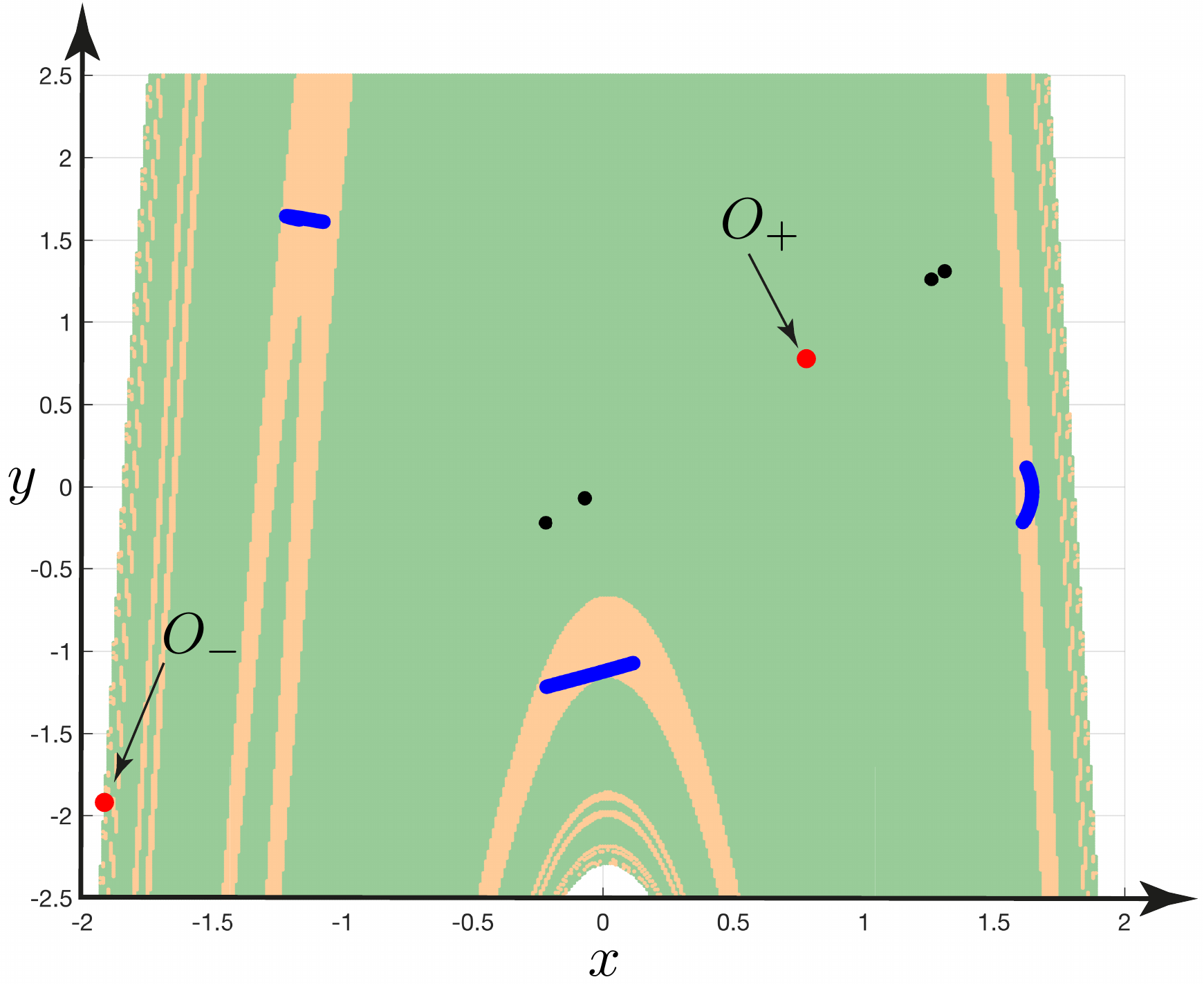}}
    \caption{Self-excited periodic attractor (black) with respect to the unstable equilibria $O_{\pm}$ (red) and its basin of attraction (green); self-excited chaotic attractor (blue) with respect to the unstable equilibrium $O_{-}$ (red), its basin of attraction (orange).}
    \label{fig::henon:attr:hid}
\end{figure}

Difficulties in obtaining a guaranteed result of UPO stabilization may be connected with complex dynamics effects such as multistability and possible existence of hidden attractor \cite{DudkowskiJKKLP-2016,PrasadIJBC-2015,LeonovK-2013-IJBC,KuznetsovLMPS-2018,KuznetsovMKK-2020,KuznetsovChua-ND2022}. For example, if we consider \eqref{henon-sys} with parameter values $a = 1.49$, $b = -0.138$ given in~\cite{DudkowskiPK-Henon-2016}, then from some initial point on the unstable manifold of the saddle $O_{-}$ a self-excited attractor with respect to $O_{-}$ can be visualized and a self-excited periodic attractor can be visualized from vicinity of $O_{+}$ (see Fig.~\ref{fig::henon:attr:hid}). In this case, other modifications of RL are required to stabilize the chaotic dynamics in the system, e.g., the actor-critic method \cite{WatkinsD-QL-1992} may be effective.
This method is inherently an extension of policy gradient approaches, in which the model tries to learn (optimize) the policy while simultaneously using it to explore the environment. It nonetheless borrows from another class of methods (Q-learning methods) a value function for the estimation of selected actions instead of rewards given directly by the environment.
We show that such an expansion of the range of AI tools makes it possible to achieve maximum progress in revealing and analyzing of the complex irregular behavior of the H\'{e}non model for various admissible values of the model parameters.
This, in turn, will enrich not only the mathematical results in this field, but also improve the connection with applications.

\section*{Conclusion}

Revealing an irregular regime, including a chaotic one, and stabilizing such complex dynamics with the help of control in economic models significantly improves the forecasting of economic processes and, thereby, optimizes the complicated management decision-making system. We considered two economic concepts represented by the continuous time OLG model and the discrete time ST model reduced to the H\'{e}non map. On the example of these models, we demonstrated the effectiveness and limitations of using the combination of DE and SOMA algorithms and the Pyragas-like method (the first model) and the continuous deep Q-learning algorithm (the second model) to solve the problem of forecasting and controlling models with irregular dynamics. We found two period-5 UPOs in OLG model via DE and SOMA algorithms, performed fine-tuning of parameters for delayed control in the Pyragas method to stabilize UPOs, conducted the optimal tuning of control parameters to maximize the basin of attraction. In a H\'{e}non-like model of pricing patterns, we stabilized period-2 UPO via the continuous deep Q-learning method, by considering specific form of the reward \eqref{r-reward}, compaired the result with direct application of the Pyragas method, and revealed practical issues of this approach in the case of multistability. Proposed approach allows us not only to identify critical states and stabilize chaotic dynamics in the models, but also to make substantial progress in testing powerful artificial intelligence tools to form recommendations for adoption policy decisions.

\section*{Acknowledgments}
The work is carried out with the financial support of the Russian Science Foundation: project 22-11-00172 (sections~1), St.Petersburg State University grant Pure ID 75207094 (sections~2), the Leading Scientific Schools program: project NSh-4196.1.1 (sections~4), and the Grant of SGS No.~SP2023/050, VŠB-Technical University of Ostrava, Czech Republic (section~3).

% \bibliographystyle{elsarticle-num}

% \bibliography{henon,2021-AI-OLG,bib_economics_AI,Ivans-bib,bib_full,bib_shap} 

\end{document}